\documentclass[12pt]{article}
\usepackage{amsmath,amsthm,amsfonts,amssymb,latexsym}

\date{\today}

\newtheorem{theorem}{Theorem}[section]

\newtheorem{lemma}[theorem]{Lemma}
\newtheorem{corollary}[theorem]{Corollary}
\newtheorem{proposition}[theorem]{Proposition}

\newtheorem{question}[theorem]{Question}

\newtheorem{example}[theorem]{Example}

\newcommand{\pr}{\mathrm{pr}}

\newcommand{\w}{\omega}

\newcommand{\B}{\mathcal{B}}
\newcommand{\C}{\mathcal{C}}

\newcommand{\A}{\mathcal{A}}

\newcommand{\TT}{\mathbb{T}}
\newcommand{\LL}{\mathcal{L}}

\newcommand{\e}{\varepsilon}

\newcommand{\osc}{\mathrm{osc}}
\newcommand{\grp}{\mathrm{grp}}
\newcommand{\sgrp}{\mathrm{sgrp}}

\newcommand{\ZFC}{\mathrm{ZFC}}
\newcommand{\MA}{\mathrm{MA}}

\begin{document}

\title{A new Lindel\"of topological group}
\author{Du\v{s}an Repov\v{s} and   Lyubomyr Zdomskyy}
\maketitle
\begin{abstract}
We show that the subsemigroup of the product of $\w_1$-many circles
generated by the $L$-space  constructed by J.~Moore is again an
$L$-space. This leads to a new example of a Lindel\"of topological
group.
 The question whether all finite powers of
this group are Lindel\"of remains open.
\end{abstract}

\section{Introduction}
This paper is devoted\footnotetext{The  authors were supported by
the Slovenian Research Agency grants P1-0292-0101-04, J1-9643-0101,
and BI-UA/07-08-001. The second author would also like to thank FWF
grant P19898-N18  for support for this research. A part of the work was done in 2007
when the second author was a Post-Doctoral Fellow at 
the Weizmann Institute of Science in Israel.

\normalsize \emph{Keywords and phrases.} Lindel\"of topological
groups, $L$-space, $L$-(semi)group, Martin's Axiom.

\emph{2000 MSC.} Primary: 22A20. Secondary: 54D20,
54H11.\footnotesize} to one of the possible  approaches to the
problem posed by Arhangel'ski\u\i\ \cite{Ar80} concerning existence
of a Lindel\"of topological group with non-Lindel\"of square. This
approach is based on  the recent deep result of Moore \cite{Mo}
asserting that there exists an $L$-space in $\ZFC$. We recall that
an $L$-space is a regular hereditarily Lindel\"of nonseparable
topological space. The connection between $L$-spaces and
preservation of Lindel\"ofness by finite powers is given by the
following result, which is a corollary of \cite[Theorem
7.10]{Tobook} and its proof:

\begin{theorem} \label{zvjazok}
Suppose that  $X$ is a regular topological space with countable
tightness  and $Y$ is a non-separable subspace of $X$. If all finite
powers of $X$ are Lindel\"of, then there exists a c.c.c. poset
$\mathbb P$ and a family $\mathcal D$ of dense subsets of $\mathbb
P$ of size $|\mathcal D|=\w_1$ such that if there exists a filter
$G\subset\mathbb P$ meeting each $D\in\mathcal D$, then $Y$ has an
uncountable discrete subspace.

Consequently, if  $\MA_{\w_1}$ holds and $X$ is a regular
topological space with countable tightness containing an
$L$-subspace, then some of the finite powers of $X$ are not
Lindel\"of.
\end{theorem}

The $L$-space constructed in \cite{Mo} is a subspace of
$\Sigma_{\w_1}$, the $\Sigma$-product of $\w_1$ many circles. It is
well-known \cite{Co} that all finite powers of this $\Sigma$-product
have countable tightness.   Theorem~\ref{zvjazok} suggests the
following open question.

\begin{question} \label{rel_arch}
Let $\LL$ be the $L$-space constructed in \cite{Mo}. Is the subgroup
of $\Sigma_{\w_1}$ generated by $\LL$ a Lindel\"of group? More
generally, can $\LL$ be embedded into a Lindel\"of subgroup $G$ of
$\Sigma_{\w_1}$?
\end{question}

The $L$-space constructed in \cite{Mo} remains an $L$ space in
extensions by a wide class of forcing notions containing all c.c.c.
ones.
 Therefore
if the answer to Question~\ref{rel_arch} is positive, i.e. $\LL$ can
be embedded into a Lindel\"of  subgroup
   $G$ of $\Sigma_{\w_1}$, then Theorem~\ref{zvjazok} would imply that
   some of the finite powers of $G$
are not Lindel\"of in ZFC.

In this paper we make a step towards the solution of
Question~\ref{rel_arch}. Using the ideas of \cite{Mo}, we show in
Section~\ref{semigroups} that the subsemigroup of  $\Sigma_{\w_1}$
generated by $\LL$ is an $L$-space. Thus there exists an
$L$-semigroup with cancellation, which seems to have not been noted
elsewhere. On the other hand, the group generated by $\LL$ contains
a copy of the one-point compactification of the discrete space of
size $\w_1$, and hence is not hereditarily
Lindel\"of. 
         In Section~\ref{prob_arch} we prove
that the subgroup of the Tychonoff product of $\w_1$-many circles
 generated by the union of $\LL$ and certain
meager $\sigma$-compact subspace is Lindel\"of, which speaks for the
positive answer to Question~\ref{rel_arch}. However, this group has
uncountable tightness, and consequently it is not within the scope
of applications of Theorem~\ref{zvjazok}.

The authors were able to find only two consistent examples of a
Lindel\"of group $G$ with
 non-Lindel\"of square in the literature, see \cite{Ma87} and \cite{To93}.
Malykhin's example is constructed under $\mathrm{cof}(\mathcal
M)=\w_1$ in terms of \cite{Bl}, while Todor\v{c}evi\'c uses the
additional assumption that there exists a countably additive measure
extending the Lebesgue measure  and which is defined  on all sets of
reals. Both of these assertions  contradict Martin's Axiom. The
existence of such a group $G$ is also consistent with $\MA$: Soukup
\cite{So01} constructed a model of  $\ZFC+\MA$ which contains an
$L$-group of countable tightness (an $L$-group is a  topological
group whose underlying topological space is an $L$-space.) Therefore
Theorem~\ref{zvjazok} implies that some of the finite powers of $G$
are not Lindel\"of.

All spaces considered here are assumed to be Tychonoff.

\section{$L$-semigroups with cancellation} \label{semigroups}

We  briefly discuss Theorem~\ref{zvjazok} before passing to
$L$-semigroups.
\medskip

\noindent\textit{Proof  sketch of Theorem~\ref{zvjazok}. \ } The
direct application of \cite[Theorem 7.10]{Tobook} gives
Theorem~\ref{zvjazok} only for spaces $X$ such that \emph{all finite
powers} of $X$ have countable tightness. However for a pair $X,Y$ of
spaces satisfying the premises of Theorem~\ref{zvjazok} one can
easily construct a continuous map $f:X\to\Sigma_{\w_1}$ such that
$f(Y)$ is not separable, see, e.g., the proof of \cite[Corollary
2.3]{MilZdo??}. Since all finite powers of $f(X)$ have countable
tightness, we can apply to $f(X), f(Y)$ the same argument as in the
proof of \cite[Theorem 7.10]{Tobook} and then pull the conclusion
back to $X,Y$. This way we get Theorem~\ref{zvjazok}. \hfill $\Box$
\medskip

In the rest of this section we follow  the notations from \cite{Mo}.
 Developing the ideas of Todor\v{c}evi\'c \cite{To89},
 Moore considered  the function $\osc:\{(\alpha,\beta)\in\w_1^2:\alpha<\beta\}\to\w$
 having  strong combinatorial properties.
We shall give   more detailed definition of this function in
Example~\ref{LindPowers}. For the purposes of this section the
following fundamental result is sufficient.

\begin{theorem}\textup{(\cite[Theorem 4.3]{Mo}).} \label{fromMoore1}
For every uncountable families of pairwise disjoint sets
$\A\subset[\w_1]^k$ and $\B\subset [\w_1]^l$  and every $n\in\w$,
there exist $a\in\A$ and $b_m\in\B$, $m<n$, such that for all $i<k$,
$j<l$, and $m<n$:

$$ a<b_m, \mbox{\ \ and\ \ }  \osc(a(i),b_m(j))=\osc(a(i), b_0(j))+m.$$
\end{theorem}
 (Here $a<b$ means $\max a< \min b$.)

Let $(z_\alpha)_{\alpha<\w_1}$ be a sequence of points  on the
circle $\TT=\{z\in\mathbb C: |z|=1\}$ which are rationally
independent. (We consider $\TT$ as a subgroup of $\mathbb
C\setminus\{0\}$ with the multiplication.) Given any
$\alpha<\beta<\w_1$, set $o(\alpha,\beta)=
z_\alpha^{\osc(\alpha,\beta)+1}$.  We define $w_\beta\in\TT^{\w_1}$
by letting \smallskip

\centerline{\(\begin{array}{c}
w_\beta(\alpha)= \left\{  \begin{array}{ll} o(\alpha,\beta), & \mbox{ if } \alpha<\beta, \\
1, & \mbox{ otherwise. }
\end{array}
 \right.
\end{array}  \) }
\medskip

It was showed in \cite[Theorem~7.11]{Mo} that for every uncountable
$X\subset\w_1$ the space $\LL_X=\{w_\beta|X:\beta\in X\}$ is an
$L$-space. The methods developed in \cite{Mo} allow one to slightly
extend this result.

For  a subset $A$ of a group $G$, we denote by $\sgrp(A)$ and
$\grp(A)$ the smallest subsemigroup and  subgroup of $G$ containing
$A$, respectively. In particular, $\sgrp(\LL_X) $ stands for the
subsemigroup of $\TT^X$ generated by $\LL_X$. A \emph{semigroup with
cancellation} is a semigroup $H$ such that both of the equalities $h
h'=h h''$ and $h''h=h'h$ imply $h'=h''$, where $h,h',h''\in H$.

\begin{theorem} \label{l-grp}
For every uncountable $X\subset\w_1$ the subspace $\sgrp(\LL_X)$ of
$\TT^X$ is an $L$-space. In particular, $\sgrp(\LL_X)$ is an
$L$-subsemigroup of $\TT^X$ with cancellation.
\end{theorem}


  The following classical result independently proved by Kronecker
and Tchebychef  will be useful.

\begin{theorem} \label{kro}
Suppose that $z_i$, $i<k$, are elements of $\TT$ which are
rationally independent. For every $\e>0$ there exists a natural
number $n_\e$ such that if $u, v\in\TT^k$, then there is an $m<n_\e$
such that $|u_i z_i^{m}-v_i|<\e$ for all $i<k$.
\end{theorem}

The next proposition resembles \cite[Theorem~5.6]{Mo}.

\begin{proposition} \label{summa}
Let $\A\subset [\w_1]^k$ and $\B\subset [\w_1]^l$ be uncountable
families of pairwise disjoint sets. Then for every sequence
$(U_i)_{i<k}$ of open subsets of $\TT$,  every partitions
$k=u_0\sqcup u_1$ and $l=v_0\sqcup v_1$, and arbitrary sequence
$(n_j)_{j<l}$ of integers with the property $\sum_{j\in v_r}n_j\neq
0$ for all $r\in 2$, there exist $a\in \A$ and $b\in\B$ such that
$a<b$ and
 $$\prod_{j\in v_r} o(a(i), b(j))^{n_j}\in U_i$$
  for all $i\in u_r$ and $r\in 2$.
\end{proposition}
\begin{proof}
There is no loss of generality in assuming that $U_i$ is an
$\e$-ball around a point $t_i$ for some fixed $\e>0$. Set
$N_r=\sum_{j\in v_r} n_j$, $r\in 2$, and let
$\delta=\e/\max\{|N_0|,|N_1|\}$.
 Passing to
an uncountable subset of $\A$, if necessary, we may additionally
assume that the numbers $n_{\delta}$ given by Theorem~\ref{kro} for
the sequence $z_{a(i)}$, $i\in k$, are the same for all $a\in\A$.

Let $a\in A$ and $b_m\in\B$, $m<n_{\delta}$, be such as in
Theorem~\ref{fromMoore1}, i.e. for all $i<k$, $j<l$, and $m<N$ we
have $a< b_m$  and $ \osc(a(i),b_m(j))=\osc (a(i), b_0(j))+m. $ For
each $r\in 2$ and $i\in u_r$ put $t'_i=\prod_{j\in v_r} o(a(i),
b_0(j))^{n_j}$. Let $t''_i$ be such that the $N_r$-th power of
$t''_i$ equals $t_i {t'_i}^{-1}$, and let $W_i$ be the $\delta$-ball
around $t''_i$, where $i\in k$. By the definition of $n_\delta$,
there exists $m<n_\delta$ such that
$$ z_{a(i)}^m\in W_i $$ for all $i<k$. Set $b=b_m$.
Then
$$ \prod_{j\in v_r} o(a(i), b(j))^{n_j}\in \prod_{j\in v_r} o(a(i),
b_0(j))^{N_r} z_{a(i)}^{m\cdot N_r}\subset t'_i W_i^{N_r}. $$ The
$W_i$'s were chosen in such a way that $W_i^{N_r}$ is a subset of
the $\e$-ball around $t_i {t'_i}^{-1}$. This completes our proof.
\end{proof}

The following proposition is reminiscent of \cite[Theorem~7.10]{Mo}.

\begin{proposition} \label{like7.10}
If $X,Y\subset\w_1$ are disjoint, then there is no continuous
injection of any uncountable subspace of $\sgrp(\LL_X)$ into
$\LL_Y$.
\end{proposition}
\begin{proof}
Suppose to the contrary  that such an injection $g$ of an
uncountable subset $Q$ of $\sgrp(\LL_X)$ into $\LL_Y$ exists.
Passing to an uncountable subset of $Q$, if necessary, we may assume
that there exist $m\in\w$, a $\Delta$-system $\C$ of subsets of $X$
of size $m$ with a root $d$, and a sequence $(n'_j)_{j<m}$ of
positive integers, such that $s_c=\prod_{j\in m} w_{c(j)}^{n'_j}\in
Q$ and $g:s_c\mapsto w_{f(c)}$, where $f:\C\to Y$ is an injection.
It is also clear that there is no loss of generality in assuming
that $d=\emptyset$.

For each $\xi<\w_1$, let $c_\xi\in \C$ and $\zeta_\xi\in Y$ be such
that $f(c_\xi)>\zeta_\xi$ and if $\xi<\xi'$, then
$c_\xi<\zeta_{\xi'}$. Let $\Theta\subset\w_1$ be uncountable such
that for some open neighborhood $V\subset\TT$,
$w_{f(c_\xi)}(\zeta_\xi)\not\in\bar{V}$ for all $\xi\in\Theta$.

Applying the continuity of $g$ at $s_c$ to $
W_\xi=\{w\in\LL_Y:w(\zeta_\xi)\not\in\bar{V}\},$
 we can find a basic open neighborhood $U_\xi$ of $s_{c_\xi}$ in
$Q$ such that $ g(U_\xi)\subset W_\xi$. Applying the $\Delta$-system
 lemma
\cite[Theorem~1.6]{Ku80}
  and the second countability of $\TT$ again, we see that
there exist $k_0\in\w$, an uncountable $\Theta'\subset\Theta$, open
neighborhoods $(U'_i)_{i\in k_0}$ in $\TT$, and $a_\xi\in [X]^{k_0}$
such that for all $\xi\in\Theta'$:

$(i)$ \ $\{a_\xi:\xi\in\Theta'\}$ is a $\Delta$-system  with a root
$a$;

$(ii)$ \ the set $\{w\in Q:\forall i<k_0 \ (w(a_\xi(i))\in U'_i)\}$
contains $s_{c_\xi}$ and is a subset of $U_\xi$;

$(iii)$ \  $|c_\xi\cap f(c_\xi)|$ does not depend on $\xi$; and

$(iv)$ \ $|(a_\xi\setminus a)\cap \zeta_\xi|$ does not depend on
$\xi$.
\medskip

Let $\A$ be the collection of all $a_\xi\cup\{\zeta_\xi\}\setminus
a$, $\xi\in\Theta'$, and let $k$ be the size of elements of $\A$.
Let also $\B$ be the collection of all $c_\xi\cup\{f(c_\xi)\}$,
where $\xi\in\Theta'$, and $l=m+1$.

Now, let $k=u_0\sqcup u_1$ and $l=v_0\sqcup v_1$
 be the partitions of $k$ and $l$ defined as
follows: $u_1=\{|(a_\xi\setminus a)\cap \zeta_\xi|\}$,
$v_1=\{|c_\xi\cap f(c_\xi)|\}$,  $u_0=k\setminus u_1$, and
$v_0=k\setminus v_1$ (conditions $(iv)$ and $(iii)$ mean that the
partitions do not depend on a particular $\xi\in\Theta'$.)
 For every $j\in l$ we put
 \smallskip

\centerline{\(\begin{array}{c}
n_j= \left\{  \begin{array}{ll} n'_j, & \mbox{ if }j<|c_\xi\cap f(c_\xi)|, \\
1, & \mbox{ if }j=|c_\xi\cap f(c_\xi)|, \\
n'_{j-1}, & \mbox{ if }j>|c_\xi\cap f(c_\xi)|.
\end{array}
 \right.
\end{array}  \) }
\medskip

Finally, for every $i\in k$ we define $U_i$ as follows:

\centerline{\(\begin{array}{c}
U_i= \left\{  \begin{array}{ll} U'_{i+|a|}, & \mbox{ if }i<|(a_\xi\setminus a)\cap \zeta_\xi|, \\
V, & \mbox{ if }i=|(a_\xi\setminus a)\cap \zeta_\xi|, \\
U'_{i+|a|-1}, & \mbox{ if }i>|(a_\xi\setminus a)\cap \zeta_\xi|.
\end{array}
 \right.
\end{array}  \) }
\medskip

Applying Proposition~\ref{summa}, it is possible to find
$\xi<\xi'\in\Theta'$ such that

\centerline{ $a=a_\xi\cup\{\zeta_\xi\}<
c_{\xi'}\cup\{f(c_{\xi'})\}=b$ and }
\smallskip

\centerline{ $\prod_{j\in v_r} o(a(i), b(j))^{n_j}\in U_i$ for all
$i\in u_r$ and $r\in 2$.}
\smallskip

The $U_i$'s and $n_j$'s were defined in such a way that the second
condition under $r=1$ gives $ w_{f(c_{\xi'})} (\zeta_\xi) \in V, $
and for $r=0$ this gives

$$ s_{c_{\xi'}}(a_{\xi}(i))=\prod_{j\in m} w_{c_{\xi'}(j)} (a_\xi (i))^{n'_j}\in U'_i $$
for all $i\geq |a|$, while for $i<|a|$ the above trivially holds by
$(i)$ and $(ii)$. But now $s_{c_{\xi'}}\in U_\xi$ even though
$g(s_{c_{\xi'}})(\zeta_\xi)= w_{f(c_{\xi'})}(\zeta_\xi)\in V$,
contradicting the choice of $U_\xi$. The proof is thus finished.
\end{proof}
\medskip

\noindent\textit{Proof of Theorem~\ref{l-grp}.} \ The ``+1'' in the
definition of the function $o$ clearly ensures that the closure in
$\sgrp(\LL_X)$ of any countable subset of $\sgrp(\LL_X)$ is
countable. Indeed, suppose that $H$ is a countable subset of $\LL_X$
and $\alpha\in\w_1$ is such that $\alpha>\xi$ for all $\xi$ with
$w_\xi|X\in H$. Thus $t(\gamma)=1$ for every $t\in\sgrp(H)$ and
$\gamma\geq\alpha$.
 Let us fix $s=\prod_{i\leq n}w_{\xi_i}^{m_i}|X\in \sgrp(\LL_X)$.
Without loss of generality, $\xi_0<\xi_1<\ldots\xi_n$ and $m_n\neq
0$. If $\xi_n>\alpha$,
$$ s(\max\{\alpha,\xi_{n-1}\})=z_{\max\{\alpha,\xi_{n-1}\}}^{m_{n}(\mathrm{osc}(\max\{\alpha,\xi_{n-1}\},\xi_n)+1)}\neq 1,$$
and consequenty $s$ is not in the closure of $\sgrp(H)$.

Therefore, if $\sgrp(\LL_X)$ were not hereditarily Lindel\"of, it
would contain an uncountable discrete subspace $Q$. The above means
that for every $q\in Q$ there exists a basic open subset $U_q\ni q$
of $\TT^{\w_1}$ such that $U_q\cap Q=\{q\}$. Since each $U_q$
depends on finitely many coordinates, we can find an uncountable
$Y\subset X$ such that $|\w_1\setminus Y|=\w_1$ and $Q|Y=\{q|Y:q\in
Q\}$ is still discrete. Then any injection
$g:Q|Y\to\LL_{\w_1\setminus Y}$ is continuous, which contradicts
Proposition~\ref{like7.10}. \hfill  $\Box$
\medskip

The following technical statement will be crucial in the next
section.

\begin{corollary} \label{tech}
Let $\C\subset[\w_1]^{l}$ be an uncountable family of pairwise
disjoint sets  and $(n_j)_{j<l}$ be a sequence of integers with
$\Sigma_{j<l}n_j\neq 0$. Then for every $X\subset\w_1$ such that
$\bigcup \C\subset X$, the subspace
$$ \{\prod_{j<l}w_{c(j)}^{n_j}|X\: :\: c\in\C\}$$
of $\Sigma_X$ is hereditarily Lindel\"of.
\end{corollary}
\begin{proof}
Almost literal repetition of the  proof of
Proposition~\ref{like7.10} (just a couple of the first lines should
be omitted) gives us that there is no continuous injection from any
uncountable subspace of $\{\prod_{j<l}w_{c(j)}^{n_j}|X\: :\:
c\in\C\}$
 into $\LL_Y$ provided $Y\cap\bigcup C=\emptyset$.
Now it suffices to apply the same argument as in the proof of
Theorem~\ref{l-grp}.
\end{proof}

In the same way we can also prove the following proposition, which
shows that it is essential in Theorem~\ref{zvjazok} to consider
finite powers and not just finite products.

\begin{proposition}
For every finite family $\{X_0,\ldots,X_n\}$ of uncountable pairwise
disjoint subsets of $\w_1$, the product
$\sgrp(\LL_{X_0})\times\cdots\times\sgrp(\LL_{X_n})$ is an
$L$-space.
\end{proposition}

On the other hand, it is easy to prove that $\grp(\LL_X)$ is not
hereditarily Lindel\"of. We shall use the following consequence of
\cite[Proposition~7.13]{Mo}.
\begin{proposition} \label{7.13}
For every $\beta<\w_1$ the set $\{w_\xi|\beta:\xi<\w_1\}$ is
countable.
\end{proposition}

For a cardinality $\tau$ we denote by $A(\tau)$ the one-point
compactification of the discrete space of size $\tau$. The following
proposition corresponds to \cite[Theorem~7.2]{Mo}.

\begin{proposition}\label{notLind}
$\grp(\LL_X)$ contains  a copy of $A(\w_1)$.
\end{proposition}
\begin{proof}
Using Proposition~\ref{7.13} we can construct  two increasing
transfinite sequences $(\xi_\beta)_{\beta<\w_1}$ and
$(\zeta_\beta)_{\beta<\w_1}$ of ordinals with the following
properties:

$(i)$ \ $\zeta_\beta>\xi_\beta$ for all $\beta<\w_1$;

$(ii)$ \ $\xi_{\beta'}>\zeta_\beta$ for all $\beta<\beta'<\w_1$; and

 $(iii)$ \
 $w_{\xi_{\beta'}}|\sup\{\zeta_\beta:\beta<\beta'\}=
w_{\zeta_{\beta'}}|\sup\{\zeta_\beta:\beta<\beta'\} $ for all
$\beta'<\w_1$.

\noindent A direct verification shows that $\{w_{\zeta_\beta}\cdot
w_{\xi_\beta}^{-1}:\beta<\w_1\}\cup\{0\}$ is a copy of $A(\w_1)$.
\end{proof}

\section{An example of a  Lindel\"of group} \label{prob_arch}

In this section we shall construct an example of a Lindel\"of group
$G$ containing $\LL_X$ of the form $\grp(\LL_X\cup K)$ for a meager
$\sigma$-compact subgroup $K$ of $\mathbb T^{\w_1}$ defined below.
This group has uncountable tightness, and hence
Theorem~\ref{zvjazok} cannot be used here to deduce that $G^n$ is
not Lindel\"of for some $n\in\w$. We do not know whether  all finite
powers of  the group $G$ constructed in Example~\ref{LindPowers} are
Lindel\"of.

Let
$$ K=\{(z_{\alpha}^{p_\alpha})_{\alpha<\w_1}:\forall \alpha<\w_1\ (p_\alpha\in\mathbb Z) \wedge (\sup\{|p_\alpha|:\alpha<\w_1\}<\infty) \}. $$
It is clear that $K$ is a meager $\sigma$-compact subgroup of
$\mathbb T^{\w_1}$. In addition, \cite[Theorem~7.14]{Mo} implies
that $\LL_X \cap \pr_{X}K$ is at most countable for every $X\in
[\w_1]^{\w_1}$.

\begin{example} \label{LindPowers}
Let $X$ be an uncountable subset of  $\w_1$ and  $G=\grp(\LL_X)
\cdot \pr_{X}K $. Then $G$ is Lindel\"of.
\end{example}

First we shall prove some auxiliary statements. At this point we
need to go a bit deeper  into the construction of the function $o$,
see \cite[Section~2]{Mo}. Summarizing Facts 1 and 2 from \cite{Mo}
we conclude that there exists a function
$L:\{(\alpha,\beta)\in\w_1^2: \alpha\leq \beta\}\to [\w_1]^{<\w}$
with the following properties:
\begin{itemize}
\item[$(i)$] $L(\alpha,\beta)\subset\alpha$ and $L(\alpha,\beta)=\emptyset$ if
and only if $\alpha=0$ or $\alpha=\beta$;
\item[$(ii)$] If $\alpha\leq\beta\leq\gamma$ and
$L(\beta,\gamma)<L(\alpha,\beta)$, then
$L(\alpha,\gamma)=L(\beta,\gamma)\cup L(\alpha,\beta)$; and
\item[$(iii)$]
If  $\beta$ is limit, then  $\lim_{\alpha\to\beta}\min
L(\alpha,\beta)=\beta$.
\end{itemize}

The definition of $o$ also involves such a standard object as a
\emph{coherent sequence} of functions, i.e. a sequence
$(e_\alpha)_{\alpha\in\w_1}$ such that $e_\alpha:\alpha\to\w$,
$e_\alpha$ is finite-to-one, and for arbitrary $\alpha<\beta$, the
set $\{\xi<\alpha:e_\alpha(\xi)\neq e_\beta(\xi)\}$ is finite. Now,
$\mathrm{osc}(\alpha, \beta)$ is the cardinality of the set
$\mathrm{Osc}(e_\alpha,e_\beta, L(\alpha,\beta))$ defined as
follows: $$\{\xi\in L(\alpha,\beta)\setminus\min
L(\alpha,\beta):e_\alpha(\xi^{-})=e_\beta(\xi^{-})\wedge
e_\alpha(\xi)>e_\beta(\xi)\}, $$ where $\xi^{-}$ is the greatest
element of $L(\alpha,\beta)$ smaller than $\xi$.

\begin{lemma} \label{l1_12.08}
Let $a\in [\w_1]^k$ and $(n_i)_{i\in k}$ be a finite sequence of
integers with the property $\sum_{i\in k} n_i=0$. Then the set $
\{\sum_{i\in k} \osc(\alpha, a(i))\cdot n_i : \alpha< a(0)\} $ is
finite.
\end{lemma}
\begin{proof}
Assuming the converse, we can find an ordinal $\eta\leq a(0)$ and a
sequence $(\xi_n)_{n\in\w}$ of ordinals converging to $\eta$ such
that $\xi_n<\xi_{n+1}$ and $|\sum_{i\in k} \osc(\xi_n, a(i))\cdot
n_i|\geq n$. Let $\gamma_0, \gamma_1<\eta$ be such that
$L(\eta,a(i))<\gamma_0$ for all $i\in k$ and
$L(\gamma,\eta)>\gamma_0$ for all $\gamma_1\leq\gamma<\eta$, and
$e_{a(i)}|(\gamma_0,\eta)=e_{a(j)}|(\gamma_0,\eta)$ for all $i,j\in
k$ (this can be done by the facts above). Then for every $i\in k$
and $\gamma_1\leq\gamma<\eta$, $L(\gamma, a(i))=L(\gamma,\eta)\cup
L(\eta, a(i))$, and hence $\mathrm{Osc}(e_\gamma, e_{a(i)},
L(\gamma, a(i))) = \mathrm{Osc}(e_\gamma, e_{a(i)},  L(\eta,
a(i))\cup L(\gamma,\eta))$. Let $q_\gamma=\mathrm{Osc}(e_\gamma,
e_{a(i)}, L(\gamma,\eta))$ (it does not depend on $i$ by our choice
of $\gamma_1$). Therefore
$$ |\mathrm{Osc}(e_\gamma, e_{a(i)}, L(\gamma, a(i)))| =
|\mathrm{Osc}(e_\gamma, e_{a(i)},  L(\eta, a(i)))| + q_\gamma +
s_\gamma, $$ where $s_\gamma\in\{0,1\}$ is the number indicating
whether $\min L(\gamma,\eta)$ is included into
$$\mathrm{Osc}(e_\gamma, e_{a(i)}, L(\gamma, a(i)))=
\mathrm{Osc}(e_\gamma, e_{a(i)},  L(\eta, a(i))\cup
L(\gamma,\eta))$$
 or not. Set  $M=\max_{i\in k}|L(\eta, a(i))|$.
Then for every $\gamma\in (\gamma_1, \eta)$ we have
\begin{eqnarray*}
|\sum_{i\in k} \osc(\gamma, a(i))\cdot n_i|= |\sum_{i\in k}
|\mathrm{Osc}(e_\gamma, e_{a(i)}, L(\gamma, a(i)))|\cdot n_i| = \\
=|\sum_{i\in k} (|\mathrm{Osc}(e_\gamma, e_{a(i)},  L(\eta,
a(i)))| + q_\gamma + s_\gamma)\cdot n_i | =\\
=|\sum_{i\in k} |\mathrm{Osc}(e_\gamma, e_{a(i)},  L(\eta,
a(i)))|\cdot n_i + \sum_{i\in k}q_\gamma\cdot n_i + \sum_{i\in k}
s_\gamma\cdot n_i | = \\
=|\sum_{i\in k} |\mathrm{Osc}(e_\gamma, e_{a(i)},  L(\eta,
a(i)))|\cdot n_i  + \sum_{i\in k} s_\gamma\cdot n_i | \leq
(kM+1)\sum_{i\in k}|n_i|,
\end{eqnarray*}
which is a contradiction.
\end{proof}

\noindent\textit{Proof of Example~\ref{LindPowers}.} Assuming that
$G$ is not Lindel\"of, we can find an increasing family
$\{U_\alpha:\alpha<\w_1\}$ of open subsets of $\mathbb T^{X}$
covering $G$ and an element $g_\alpha\in G\setminus U_\alpha$. Using
the standard $\Delta$-system argument, we can find an uncountable
family $B\subset [X]^l$  of pairwise disjoint sets, a sequence
$(n_j)_{j<l}$ of integers, $x\in\grp(\LL_X)$, and $\{y_b:b\in
B\}\subset \pr_{X}K$ such that
\begin{equation*} x \cdot\{
\prod_{j\in l}w_{b(j)}^{n_j}\cdot y_b:b\in B\} \subset
\{g_\alpha:\alpha<\w_1\}, \label{unc1}
\end{equation*}
and hence\footnote{Formally, we should have written $w_{-}|X$
instead of $w_{-}$ throughout the proof.} the intersection
\begin{equation*}
\big(x\cdot\{ \prod_{j\in l}w_{b(j)}^{n_j}:b\in B\}\cdot
\pr_{X}K\big) \cap \{g_\alpha:\alpha<\w_1\} \label{unc2}
\end{equation*}
 is uncountable. Two cases are possible:

Case 1. $\sum_{j\in l}n_j\neq 0$. In this case  Corollary~\ref{tech}
implies that $\{ \prod_{j\in l}w_{b(j)}^{n_j}:b\in B\}$ is
hereditarily Lindel\"of,  hence $x\cdot\{ \prod_{j\in
l}w_{b(j)}^{n_j}:b\in B\}\cdot \pr_{X}K$ is Lindel\"of being a
continuous image of a product of a Lindel\"of  space with a
$\sigma$-compact, and therefore
 this set is contained in some $U_\xi$, which contradicts
the fact that it contains uncountably many $g_\alpha$'s.

Case 2. $\sum_{j\in l}n_j = 0$. Passing to an uncountable subset of
$B$, if necessary, we can additionally assume that
$B=\{b_\xi:\xi<\w_1\}$ and $ b_\xi> b_\eta$ provided that
$\eta<\xi$. Let $y'_b(\alpha)=\prod_{j\in l}w_{b(j)}^{n_j}(\alpha)$
if $\alpha<b(0)$, and $y'_b(\alpha)=1$ otherwise, where $b\in B$.
Applying Lemma~\ref{l1_12.08}, we conclude that $y'_b\in K$ for all
$b\in B$. In addition, it is easy to see that $C=\{\prod_{j\in
l}w_{b_\xi(j)}^{n_j}\cdot (y'_{b_\xi})^{-1}:\xi<\w_1\}\bigcup\{1\}$
is a copy of $A(\w_1)$. Therefore \begin{eqnarray*} x\cdot\{
\prod_{j\in l}w_{b(j)}^{n_j}:b\in B\}\cdot \pr_{X}K=\\= x\cdot\{
\prod_{j\in l}w_{b(j)}^{n_j}\cdot (y'_b)^{-1} \cdot y'_b:b\in
B\}\cdot \pr_{X}K \subset x\cdot C \cdot \pr_{X}K\cdot \pr_{X}K,
\end{eqnarray*}
and the latter set is  a $\sigma$-compact subset of $G$, and hence
it is contained in some $U_\alpha$, which is a contradiction. \hfill
$\Box$
\medskip

\noindent\textbf{Acknowledgement.} The authors would like to thank
Boaz Tsaban and the Hungarian topological team for for many
stimulating conversation during the 1st European Set Theory Meeting
in Poland in 2007. We are particularly grateful to Justin Moore  for
the discussions regarding Theorem~\ref{zvjazok}.

Du\v{s}an Repov\v{s}, Faculty of  Mathematics and Physics and
Faculty of Education,  University of Ljubljana, Jadranska 19,
Ljubljana, Slovenija 1000.
\smallskip

\noindent \textit{E-mail address:}   \texttt{dusan.repovs@guest.arnes.si}\\
\noindent \textit{URL:} \texttt{http://www.fmf.uni-lj.si/$\tilde{\
}$repovs/index.htm}

\medskip

 Lyubomyr Zdomskyy, Kurt G\"odel Research Center for Mathematical Logic,
University of Vienna,  W\"ahringer Stra\ss e 25,  A-1090 Wien,
Austria.
\smallskip

\noindent \textit{E-mail address:}   \texttt{lzdomsky@gmail.com}\\
\noindent \textit{URL:}
\texttt{http://www.logic.univie.ac.at/$\tilde{\ }$lzdomsky}

\end{document}